\documentclass[11pt, a4paper]{article}
\usepackage{authblk}
\usepackage[utf8]{inputenc}
\usepackage{amsmath}
\usepackage{amssymb}
\usepackage{mathtools}
\usepackage{amsthm}

\newtheorem{lemma}{Lemma}[section]

\begin{document}

\title{A Second Raabe's Test and Other Series Tests}
\markright{Second Raabe's Test}
\author{Edward Huynh}
\affil{Department of Mathematical Sciences, University of Nevada, Las Vegas, Las Vegas, Box 454020, NV, USA}

\maketitle

\begin{abstract}
The classical D'Alembert's Ratio Test is a powerful test that we learn from calculus to determine convergence for a series of positive terms. Its range of applicability and ease of computation makes this test extremely appealing. However, it admits an inconclusive case when the limiting ratio of the terms equals 1. Several series tests like Raabe's and Gauss' Tests have been proposed in order to address this case. These tests were also generalized by Kummer through Kummer's Test. More recently, a Second Ratio Test was constructed that similarly possessed an inconclusive case. This article will present a survey of existing series tests. Secondly, it will introduce an extension of Raabe's Test to the Second Ratio Test. Thirdly, other extensions of classical tests such as Gauss' Test and Kummer's Test are proposed. Finally, it will also present proofs for the aforementioned tests and a brief application of the Second Raabe's Test.
\end{abstract}

\noindent
Infinite series exist ubiquitously in science and engineering due to their significance in applications. In the 18th Century, French Mathematician Jean Le Rond D’Alembert formulated D’Alembert’s Ratio Test, otherwise known today as the Ratio Test. D’Alembert felt that the theory of limits needed to be established rigorously and defined the notion of derivative by using the “limit of a quotient of increments” \cite{d1768reflexions}. His work on limits led to the development of the Ratio Test. The Ratio Test continues to be taught in the present in college calculus classes when studying the theory of infinite series. Even though the Ratio Test is a powerful tool in determining convergence of series, this test fails to determine convergence when the limit of the ratio is 1. This inconclusive case led to the development of several special convergence tests. German mathematician Joseph Ludwig Raabe published in 1834 a convergence test known today as Raabe’s Test \cite{NotezurTheoriederConvergenzundDivergenzderReihen}. It served as an extension for D’Alembert’s Ratio Test. Meanwhile, another German mathematician Carl Friedrich Gauss published his work on the hypergeometric series as well as a convergence test that would be known today as Gauss’ Test in 1812. This test, unlike D’Alembert’s Ratio Test and Raabe’s Test, does not have an inconclusive case due to a slightly weaker assumption on the expression of the ratios.

From calculus courses, the first tool for testing convergence of series is D’Alembert’s Ratio Test which states that, for series $\sum a_n$ of nonnegative numbers, the series converges if $\lim_{n\to \infty}\frac{a_{n+1}}{a_n} < 1$ and diverges if $\lim_{n\to \infty}\frac{a_{n+1}}{a_n} >1$. The test is inconclusive if $\lim_{n\to \infty}\frac{a_{n+1}}{a_n} = 1$. This test provides an easy and direct method for determining convergences of series. However, the inconclusive case requires more sophisticated tests such as Raabe’s Test and Gauss’ Test, which are specifically designed to address the inconclusive case of the Ratio Test. 

In 2008, Ali \cite{Ali2008} generalized the Ratio Test and developed a new Ratio Test, called the Second Ratio Test, which uses two ratios, i.e. $\frac{a_{2n}}{a_n}$  and $\frac{a_{2n+1}}{a_n}$. However, the Second Ratio Test may also be inconclusive when $\frac{1}{2}$ lies between the lower limits of both ratios and the upper limits of both ratios. Along the direction of Ali’s work, one of the goals of this paper is to develop new tests to deal with the inconclusive case of the Second Ratio Test. Secondly, it will establish several other tests related to the Second Ratio Test that can be thought of as analogues of several mainstream tests that address the inconclusive case of D'Alembert's Ratio Test.

This paper will be organized in the following manner: it will begin with a brief survey of classical and recent series convergence tests; the second part of the paper will introduce the main result(s) in this paper along with their corresponding proofs; thirdly, it will cover a brief example that demonstrates an application of the main result; finally, a brief conclusion will follow that analyzes the main result and any potential research directions.

\section{Brief Survey.}

We survey the series convergence tests that have arisen in the literature. We begin with stating the classical ratio test found in standard calculus textbooks:

\textbf{Theorem 1.1} (D'Alembert's Ratio Test)
Let $\{a_n\}$ be a sequence of positive numbers. Then
\begin{itemize}
    \item[(i)] If $\limsup_{n\to \infty} \frac{a_{n+1}}{a_n} < 1$, then $\sum_{n=1}^\infty a_n$ converges.
    
    \item[(ii)] If $\liminf_{n\to \infty} \frac{a_{n+1}}{a_n} > 1$, then $\sum_{n=1}^\infty a_n$ diverges.
    
    \item[(iii)] If $\liminf_{n\to \infty} \frac{a_{n+1}}{a_n} \leq 1 \leq \limsup_{n\to \infty} \frac{a_{n+1}}{a_n}$, the test is inconclusive.
\end{itemize}
In conjunction with the root test, this test is often used for evaluation of the convergence of infinite series. Nevertheless, the test is not perfect since there exist series where the inconclusive case occurs, e.g. $a_n = \frac{1}{n}$. As a result, several tests were developed to address this case such as Raabe's Test.

\textbf{Theorem 1.2} (Raabe's Test)
Let $\{a_n\}$ be a sequence of positive numbers. Suppose we have
\begin{align*}
    \frac{a_{n+1}}{a_n} = 1- \frac{\beta}{n} + \frac{\epsilon(n)}{n},
\end{align*}
where $\beta$ is independent of $n$ and $\epsilon(n) \to 0$, as $n\to \infty$. Then
\begin{itemize}
    \item[(i)] $\sum_{n=1}^\infty a_n$ converges if $\beta > 1$.
    
    \item[(ii)] $\sum_{n=1}^\infty a_n$ diverges if $\beta < 1$.
    
    \item[(i)] If $\beta = 1$, then test is inconclusive.
\end{itemize}

\textbf{Remark.}
Notice that we may equivalently pose the test as: $\sum_{n=1}^\infty a_n$ converges if $\lim_{n\to \infty} n\left(1 - \frac{a_{n+1}}{a_n}\right) > 1$ and $\sum_{n=1}^\infty a_n$ diverges if $\lim_{n\to \infty} n\left(1 - \frac{a_{n+1}}{a_n}\right) < 1$. When $\lim_{n\to \infty} n\left(1 - \frac{a_{n+1}}{a_n}\right) = 1$, this corresponds to the inconclusive case.

By weakening the assumption on the sequence term in Raabe's Test, one can prove a similar test due to Gauss.

\textbf{Theorem 1.3} (Gauss' Test)
Let $\{a_n\}$ be a sequence of positive numbers. Suppose we have
\begin{align*}
    \frac{a_{n+1}}{a_n} = 1- \frac{\beta}{n} + \frac{\gamma(n)}{n^{1+\lambda}},
\end{align*}
where $\lambda > 0$, $\beta$ is independent of $n$ and $\gamma(n)$ is bounded. Then
\begin{itemize}
    \item[(i)] $\sum_{n=1}^\infty a_n$ converges if $\beta > 1$.
    
    \item[(ii)] $\sum_{n=1}^\infty a_n$ diverges if $\beta \leq 1$.
\end{itemize}

We note that this test sufficiently determines convergence for any positive term series provided one can formulate the ratio term as in the expression given for Gauss' Test. 

There exists a test due to Kummer which contains as special cases the Ratio and Raabe's Test:

\textbf{Theorem 1.4} (Kummer's Test)
Let $\{a_n\}$ be a sequence of positive numbers.
\begin{itemize}
    \item[(i)] Suppose there exists a sequence of positive numbers $\{p_n\}$ and a positive number $r$ such that
    \begin{align*}
        p_n\frac{a_n}{a_{n+1}} - p_{n+1} > r
    \end{align*}
    holds for sufficiently large $n$, then $\sum_{n=1}^\infty a_n$ converges.
    
    \item[(ii)] On the other hand, if there exists a sequence of positive numbers $\{p_n\}$ where $\sum_{n=1}^\infty \frac{1}{p_n}$ diverges, such that
    \begin{align*}
        p_n\frac{a_n}{a_{n+1}} - p_{n+1} < r
    \end{align*}
    for sufficiently large $n$, then $\sum_{n=1}^\infty a_n$ diverges.
\end{itemize}

By letting $p_n = 1$ or $n$, the test yields D'Alembert's Ratio Test or Raabe's Test, respectively. This test allows more flexibility in developing a series test by selecting appropriate sequences for $p_n$. While the statement of the test is sufficient for convergence and divergence, Tong \cite{Tong1994} proves that the Kummer's Test in fact provides a characterization of convergence and divergence of series.

Ali \cite{Ali2008} constructed the Second Ratio Test which can be used as a substitute for Raabe's and Gauss' Tests:

\textbf{Theorem 1.5} (Second Ratio Test)
Let $\{a_n\}$ be a sequence of positive terms. Let
\begin{align*}
    L = \max \{\limsup_{n\to \infty} \frac{a_{2n}}{a_n}, \limsup_{n\to \infty} \frac{a_{2n+1}}{a_n}\}
\end{align*}
and 
\begin{align*}
    l = \min \{\liminf_{n\to \infty} \frac{a_{2n}}{a_n}, \liminf_{n\to \infty} \frac{a_{2n+1}}{a_n}\}.
\end{align*}
Then
\begin{itemize}
    \item[(i)] If $L < \frac{1}{2}$, then $\sum_{n=1}^\infty a_n$ converges.
    
    \item[(ii)] If $l > \frac{1}{2}$, then $\sum_{n=1}^\infty a_n$ diverges.
    
    \item[(iii)] If $l\leq \frac{1}{2} \leq L$, then the test is inconclusive.
\end{itemize}

Ali shows that the Second Ratio Test may be applied to a wide range of series that may require Raabe's or Gauss' Test \cite{Ali2008}. The advantage of this test is due to its computational efficiency for certain series in comparison to these tests. 

Next, we list results which are convenient for the proofs of Second Raabe Test and other associated tests. For monotonically decreasing sequences, we have the following result by Cauchy:

\textbf{Theorem 1.6} (Cauchy's Condensation Test)
Let $\{a_n\}$ be a monotone sequence of positive numbers.

Then $\sum_{n=1}^\infty a_n$ converges if and only if $\sum_{n=1}^\infty 2^na_{2^n}$ converges.
Next, Ali \cite{Ali2008} proves a version of the Ratio Comparison Test using the two ratios in the Second Ratio Test:

\textbf{Theorem 1.7} (Second Ratio Comparison Test)
Let $\{a_n\}$ and $\{b_n\}$ be sequences of positive numbers. Suppose
\begin{align*}
    \frac{a_{2n}}{a_n} \leq \frac{b_{2n}}{b_n} \quad \text{and} \quad \frac{a_{2n+1}}{a_n} \leq \frac{b_{2n+1}}{b_n}
\end{align*}
for all large $n$. Then
\begin{itemize}
    \item[(i)] $\sum_{n=1}^\infty a_n$ converges if $\sum_{n=1}^\infty b_n$ converges.
    
    \item[(ii)] $\sum_{n=1}^\infty b_n$ diverges if $\sum_{n=1}^\infty a_n$ diverges.
\end{itemize}

\section{Main Results.}
\textbf{Theorem 2.1} (Second Raabe's Test)
Suppose $\{a_n\}$ is a sequence of positive numbers and $\lambda  \neq 0$. Let
\begin{align*}
    \liminf_{n\to \infty} \ln{\lambda n}\left( \frac{1}{2} - \frac{a_{2n}}{a_n}\right) = m_1,\quad \liminf_{n\to \infty} \ln{\lambda n}\left( \frac{1}{2} - \frac{a_{2n+1}}{a_n}\right) = m_2,\\
    \limsup_{n\to \infty} \ln{\lambda n}\left( \frac{1}{2} - \frac{a_{2n}}{a_n}\right) = M_1,\quad \limsup_{n\to \infty} \ln{\lambda n}\left( \frac{1}{2} - \frac{a_{2n+1}}{a_n}\right) = M_2.
\end{align*}
Define $m = \min \{m_1, m_2\}$ and $M= \max \{M_1, M_2\}$.
\begin{itemize}
    \item[(i)] If $m > \frac{\ln{2}}{2}$, then $\sum_{n=1}^\infty a_n$ converges.
    
    \item[(ii)] If $M < \frac{\ln{2}}{2}$, then $\sum_{n=1}^\infty a_n$ diverges.
    
    \item[(iii)] If $m\leq \frac{\ln{2}}{2} \leq M$, the test is inconclusive.
\end{itemize}

\textbf{Remark.}
We note that the Second Raabe's Test involves $\ln{n}$ as opposed to $n$ in the Raabe's Test. 

\begin{proof}
Let us first prove a lemma:
\begin{lemma}
Let $p$ be any real number. Then
\begin{align} \label{eq1}
    \lim_{n\to \infty} \ln{\lambda n}\left(\frac{1}{2}- \frac{(\ln{n})^p}{2(\ln{2n})^p}\right) = \frac{p\ln{2}}{2}
\end{align}
and
\begin{align} \label{eq2}
    \lim_{n\to \infty} \ln{\lambda n}\left(\frac{1}{2}- \frac{n(\ln{n})^p}{(2n+1)(\ln({2n+1}))^p}\right) = \frac{p\ln{2}}{2}
\end{align}
\end{lemma}
We will only prove (\ref{eq1}) since the proof for (\ref{eq2}) will follow similarly. We rewrite (\ref{eq1}) and use L'Hopital's rule:
\begin{align*}
    \lim_{n\to \infty} \ln{\lambda n}\left(\frac{1}{2}- \frac{(\ln{n})^p}{2(\ln{2n})^p}\right) &= \lim_{n\to \infty} \frac{1}{2}\left(\frac{1-\left(\frac{\ln{n}}{\ln{2n}}\right)^p}{\frac{1}{\ln{\lambda n}}}\right) \\
    &= \frac{1}{2}\lim_{n\to \infty} \left( \frac{-p \left(\frac{\ln{n}}{\ln{2n}} \right)^{p-1} \left( \frac{\frac{\ln{2n}-\ln{n}}{n}}{(\ln{2n})^2}\right)}{-\frac{1}{(\ln{\lambda n})^2} \left( \frac{1}{n} \right)} \right) \\
    &= \frac{1}{2}\lim_{n\to \infty} \left(p\left(\frac{\ln{n}}{\ln{2n}}\right)^{p-1}(\ln{2})\left(\frac{\ln{\lambda n}}{\ln{2n}}\right)^2 \right) \\
    &= \frac{p\ln{2}}{2}.
\end{align*}
The lemma has been proved.

Suppose $m$ and $M$ are defined as in the statement of the theorem. We begin by proving \textit{(i)}. Let $m > \frac{\ln{2}}{2}$. WLOG, we assume $m = m_1$. Therefore, we have $m_1 \leq m_2$. Let $\alpha, \beta \in \mathbb{R}$ such that $\frac{\ln{2}}{2} < \beta < \alpha < m$. By assumption, there exists $N_1 \geq 0$ such that $\ln{\lambda n}\left( \frac{1}{2} - \frac{a_{2n}}{a_n}\right) > \alpha$, for all $n\geq N_1$. Under algebraic rearrangement, we obtain $\frac{a_{2n}}{a_n} < \frac{1}{2} - \frac{\alpha}{\ln{\lambda n}}$, for $n \geq N_1$.

Define $b_n = \frac{1}{n\left(\ln{n} \right)^{\frac{2\beta}{\ln{2}}}}$. This allows us to obtain $\frac{b_{2n}}{b_n} =  \frac{(\ln{n})^{\frac{2\beta}{\ln{2}}}}{2(\ln{2n})^{\frac{2\beta}{\ln{2}}}}$. We invoke the lemma to obtain
\begin{align*}
    \lim_{n\to \infty} \ln{\lambda n}\left(\frac{1}{2} - \frac{b_{2n}}{b_n} \right) &= \lim_{n\to \infty}\ln{\lambda n}\left(\frac{1}{2} - \frac{(\ln{n})^{\frac{2\beta}{\ln{2}}}}{2(\ln{2n})^{\frac{2\beta}{\ln{2}}}} \right) \\
    &= \beta < \alpha.
\end{align*}
Thus, there exists $N_2 > 0$ such that $\ln{\lambda n}\left(\frac{1}{2} - \frac{b_{2n}}{b_n} \right) < \alpha$ for all $n \geq N_2$. Consequently, we have $\frac{1}{2} - \frac{\alpha}{\ln{\lambda n}} < \frac{b_{2n}}{b_n}$.

Let $N_1^* = \max \{N_1, N_2\}$. Thus, for all $n \geq N_1^*$, we obtain $\frac{a_{2n}}{a_n} < \frac{1}{2} - \frac{\alpha}{\ln{\lambda n}} < \frac{b_{2n}}{b_n}$ which implies $\frac{a_{2n}}{a_n} \leq \frac{b_{2n}}{b_n}$ for all $n\geq N_1^*$. 

We shall next show that this same inequality holds for the corresponding ratios $\frac{a_{2n+1}}{a_n}$ and $\frac{b_{2n+1}}{b_n}$. Since $m_2 > m_1 > \alpha$, by assumption there must exist $N_3 \geq 0$ such that $\ln{\lambda n}\left(\frac{1}{2} - \frac{a_{2n+1}}{a_n}\right) > \alpha$ for all $n\geq N_3$. From here, we obtain $\frac{a_{2n+1}}{a_n} < \frac{1}{2} - \frac{\alpha}{\ln{\lambda n}}$ for all $n \geq N_3$.

Using the same $b_n$ as defined earlier, we get $\frac{b_{2n+1}}{b_n} = \frac{n(\ln{n})^{\frac{2\beta}{\ln{2}}}}{(2n+1)(\ln{(2n+1)})^{\frac{2\beta}{\ln{2}}}}$. Using the lemma again, we have
\begin{align*}
    \lim_{n\to \infty} \ln{\lambda n}\left( \frac{1}{2} - \frac{b_{2n+1}}{b_n}\right) = \beta < \alpha.
\end{align*}
Thus, there exists $N_4$ such that $\ln{\lambda n}\left(\frac{1}{2} - \frac{b_{2n+1}}{b_n}\right) < \alpha$ for all $n\geq N_4$. Hence, we obtain $\frac{1}{2} - \frac{\alpha}{\ln{\lambda n}} < \frac{b_{2n+1}}{b_n}$ for all $n\geq N_4$.

If we let $N_2^* = \max \{N_3, N_4\}$, we have $\frac{a_{2n+1}}{a_n} \leq \frac{b_{2n+1}}{b_n}$ for all $n\geq N_2^*$. Finally, if we let $N = \max \{N_1^*, N_2^*\}$, we have that both $\frac{a_{2n}}{a_n} \leq \frac{b_{2n}}{b_n}$ and $\frac{a_{2n+1}}{a_n} \leq \frac{b_{2n+1}}{b_n}$ hold for $n \geq N$.

Notice that $\sum_{n=2}^\infty b_n$ converges since $\sum_{n=2}^\infty \frac{1}{n(\ln{n})^p}$ converges for any $p > 1$. Thus, by the Second Ratio Comparison Test, $\sum_{n=1}^\infty a_n$ must converge. Hence, we have proved \textit{(i)}.

We shall now prove \textit{(ii)}. The proof will more or less imitate the proof for \textit{(i)}. Suppose $M < \frac{\ln{2}}{2}$. WLOG, we assume $M = M_1$. Therefore, we have $M_2 \leq M_1$. Let $\alpha, \beta \in \mathbb{R}$ such that $M < \alpha < \beta < \frac{\ln{2}}{2}$. By assumption, there exists $N_1 \geq 0$ such that $\ln{\lambda n}\left( \frac{1}{2} - \frac{a_{2n}}{a_n}\right) < \alpha$, for all $n\geq N_1$. Under algebraic rearrangement, we obtain $\frac{1}{2} - \frac{\alpha}{\ln{\lambda n}}< \frac{a_{2n}}{a_n}$, for $n \geq N_1$.

Define $b_n = \frac{1}{n\left(\ln{n} \right)^{\frac{2\beta}{\ln{2}}}}$. This allows us to obtain $\frac{b_{2n}}{b_n} =  \frac{(\ln{n})^{\frac{2\beta}{\ln{2}}}}{2(\ln{2n})^{\frac{2\beta}{\ln{2}}}}$. We invoke the lemma to obtain
\begin{align*}
    \lim_{n\to \infty} \ln{\lambda n}\left(\frac{1}{2} - \frac{b_{2n}}{b_n} \right) &= \lim_{n\to \infty}\ln{\lambda n}\left(\frac{1}{2} - \frac{(\ln{n})^{\frac{2\beta}{\ln{2}}}}{2(\ln{2n})^{\frac{2\beta}{\ln{2}}}} \right) \\
    &= \beta > \alpha.
\end{align*}
Thus, there exists $N_2 > 0$ such that $\ln{\lambda n}\left(\frac{1}{2} - \frac{b_{2n}}{b_n} \right) > \alpha$ for all $n \geq N_2$. Consequently, we have $\frac{b_{2n}}{b_n} < \frac{1}{2} - \frac{\alpha}{\ln{\lambda n}}$.

Let $N_1^* = \max \{N_1, N_2\}$. Thus, for all $n \geq N_1^*$, we obtain $\frac{b_{2n}}{b_n} < \frac{1}{2} - \frac{\alpha}{\ln{\lambda n}} < \frac{a_{2n}}{a_n}$ which implies $\frac{b_{2n}}{b_n} \leq \frac{a_{2n}}{a_n}$ for all $n\geq N_1^*$. 

We shall next show that this same inequality holds for the corresponding ratios $\frac{a_{2n+1}}{a_n}$ and $\frac{b_{2n+1}}{b_n}$. Since $M_2 \leq M_1 < \alpha$, by assumption there must exist $N_3 \geq 0$ such that $\ln{\lambda n}\left(\frac{1}{2} - \frac{a_{2n+1}}{a_n}\right) < \alpha$ for all $n\geq N_3$. From here, we obtain $\frac{a_{2n+1}}{a_n} > \frac{1}{2} - \frac{\alpha}{\ln{\lambda n}}$ for all $n \geq N_3$.

Defining $b_n$ as earlier, we get $\frac{b_{2n+1}}{b_n} = \frac{n(\ln{n})^{\frac{2\beta}{\ln{2}}}}{(2n+1)(\ln{(2n+1)})^{\frac{2\beta}{\ln{2}}}}$. Using the lemma again, we have
\begin{align*}
    \lim_{n\to \infty} \ln{\lambda n}\left( \frac{1}{2} - \frac{b_{2n+1}}{b_n}\right) = \beta > \alpha.
\end{align*}
Thus, there exists $N_4$ such that $\ln{\lambda n}\left(\frac{1}{2} - \frac{b_{2n+1}}{b_n}\right) > \alpha$ for all $n\geq N_4$. Hence, we obtain $\frac{b_{2n+1}}{b_n}<\frac{1}{2} - \frac{\alpha}{\ln{\lambda n}}$ for all $n\geq N_4$.

If we let $N_2^* = \max \{N_3, N_4\}$, we have $\frac{b_{2n+1}}{b_n} \leq \frac{a_{2n+1}}{a_n}$ for all $n\geq N_2^*$. Finally, if we let $N = \max \{N_1^*, N_2^*\}$, we have that both $\frac{b_{2n}}{b_n} \leq \frac{a_{2n}}{a_n}$ and $\frac{b_{2n+1}}{b_n} \leq \frac{a_{2n+1}}{a_n}$ hold for $n \geq N$.

Notice that $\sum_{n=2}^\infty b_n$ diverges since $\sum_{n=2}^\infty \frac{1}{n(\ln{n})^p}$ diverges for any $p < 1$. Thus, by the Second Ratio Comparison Test, $\sum_{n=1}^\infty a_n$ diverges. Hence, we have proved \textit{(ii)}.
\end{proof}
If the limits for the ratios exist, then we may simplify Second Raabe's Test.

\textbf{Corollary 2.1.1}
Let $\{a_n\}$ be a sequence of positive numbers. Suppose
\begin{align*}
    \frac{a_{2n}}{a_n} = \frac{1}{2} - \frac{\beta_1}{\ln{\lambda n}} + \frac{\epsilon_1(n)}{\ln{\lambda n}}
\end{align*}
and
\begin{align*}
    \frac{a_{2n+1}}{a_n} = \frac{1}{2} - \frac{\beta_2}{\ln{\lambda n}} + \frac{\epsilon_2(n)}{\ln{\lambda n}}    
\end{align*}
where $\beta_1$ and $\beta_2$ are independent of $n$ and $\epsilon_1(n), \epsilon_2(n) \to 0$, as $n\to \infty$.

Define $m = \min \{\beta_1, \beta_2\}$ and $M = \max \{\beta_1, \beta_2\}$.
\begin{itemize}
    \item[(i)] If $m > \frac{\ln{2}}{2}$, then $\sum_{n=1}^\infty a_n$ converges.
    
    \item[(ii)] If $M < \frac{\ln{2}}{2}$, then $\sum_{n=1}^\infty a_n$ diverges.
\end{itemize}

Using this Corollary, we may deduce another form of the test that can be perceived as a special case of a "Second Gauss' Test".

\textbf{Corollary 2.1.2}
Let $\{a_n\}$ be a monotone decreasing sequence of positive numbers. Suppose
\begin{align*}
    \frac{a_{2n}}{a_n} = \frac{1}{2} - \frac{\beta}{\ln{n}} + \frac{\gamma(n)}{(\ln{n})^{p}},
\end{align*}
where $p >1$, $\beta$ is independent of $n$, and $\gamma(n)$ is bounded. Then
\begin{itemize}
    \item[(i)] If $\beta > \frac{\ln{2}}{2}$, then $\sum_{n=1}^\infty a_n$ converges.
    
    \item[(ii)] If $\beta \leq \frac{\ln{2}}{2}$, then $\sum_{n=1}^\infty a_n$ diverges.
\end{itemize}
\begin{proof}
Consider the series $\sum_{n=1}^\infty b_n$, where $b_n = 2^n a_{2^n}$. Then we have
\begin{align*}
    \frac{b_{n+1}}{b_n} = \frac{2a_{2(2^n)}}{a_{2^n}} &= 2\left(\frac{1}{2} - \frac{\beta}{\ln{2^n}} + \frac{\gamma(2^n)}{(\ln{2^n})^p} \right) \\
    &= 2\left(\frac{1}{2} - \frac{\beta}{n\ln{2}} +  \frac{\gamma(2^n)}{(\ln{2^n})^p} \right) \\
    &= 1 - \frac{\frac{2\beta}{\ln{2}}}{n} +  \frac{2\gamma(2^n)}{n^p\left(\ln{2}\right)^p}.
\end{align*}
By the Cauchy Condensation Test, we conclude $\sum_{n=1}^\infty a_n$ and $\sum_{n=1}^\infty b_n$ either both converge or diverge. We now consider various cases for $\beta$.

If $\beta > \frac{\ln{2}}{2}$, then $\frac{2\beta}{\ln{2}} > 1$. By Gauss' Test, $\sum_{n=1}^\infty b_n$ converges, which implies $\sum_{n=1}^\infty a_n$ converges.

On the other hand, if $\beta \leq \frac{\ln{2}}{2}$, then $\frac{2\beta}{\ln{2}} \leq 1$. Another application of Gauss's Test demonstrates that $\sum_{n=1}^\infty b_n$ diverges. This implies that $\sum_{n=1}^\infty a_n$ must also diverge.
\end{proof}

\textbf{Example 1. (The harmonic series)} Consider $a_n = \frac{1}{n}$. Notice that
\begin{align*}
    \frac{a_{2n}}{a_{n}} = \frac{1}{2} \to \frac{1}{2},\quad \text{as $n\to \infty$}
\end{align*}
and
\begin{align*}
    \frac{a_{2n+1}}{a_{n}} = \frac{n}{2n+1} \to \frac{1}{2},\quad \text{as $n\to \infty$}.
\end{align*}
Thus, the Second Ratio Test fails to give a conclusion. We will use the Second Raabe's Test:
\begin{align*}
    \lim_{n\to \infty}\ln{n}\left( \frac{1}{2} - \frac{1}{2}\right) = 0
\end{align*}
and 
\begin{align*}
    \lim_{n\to \infty} \ln{n}\left( \frac{1}{2} - \frac{n}{2n+1}\right) = 0.
\end{align*}
Since $M = 0$, we conclude by the Second Raabe's Test that $\sum a_n$ diverges.
    
Alternatively, we may use Corollary 2.1.2. where $\beta = 0$ to conclude $\sum a_n$ diverges.

We now state and prove a generalized Kummer's Test that uses the quantities involved with the Second Ratio Test:

\textbf{Theorem 2.2} (Second Kummer's Test)
Let $\{a_n\}$ be a sequence of positive numbers.
\begin{itemize}
    \item[(i)] Suppose there exists a sequence of positive numbers $\{p_n\}$ and positive numbers $r_1$ and $r_2$ such that
    \begin{align*}
        p_n\frac{a_n}{a_{2n}} - 2p_{2n} > r_1
    \end{align*}
    and
    \begin{align*}
        p_n\frac{a_n}{a_{2n+1}} - 2p_{2n+1} > r_2
    \end{align*}
    where both inequalities hold for sufficiently large $n$, then $\sum_{n=1}^\infty a_n$ converges.
    
    \item[(ii)] On the other hand, if there exists a sequence of positive numbers $\{p_n\}$ where $\sum_{k=1}^\infty \frac{1}{2^kp_{2^k}}$ diverges, and at least one of the following inequalities holds for sufficiently large $n$
    \begin{equation} \label{eq3}
        p_n\frac{a_n}{a_{2n}} - 2p_{2n} < 0
    \end{equation}
    or
    \begin{equation} \label{eq4}
        p_n\frac{a_n}{a_{2n+1}} - 2p_{2n+1} < 0
    \end{equation}
    then $\sum_{n=1}^\infty a_n$ diverges.
\end{itemize}
\begin{proof}
We begin by proving \textit{(i)}. By assumption, there exists $N_1, N_2 \in \mathbb{N}$ such that
\begin{align*}
    p_n\frac{a_n}{a_{2n}} - 2p_{2n} > r_1,\quad \forall n\geq N_1
\end{align*}
and
\begin{align*}
    p_n\frac{a_n}{a_{2n+1}} - 2p_{2n+1} > r_2, \quad \forall n\geq N_2.
\end{align*}
Define $r = \min \{r_1, r_2\}$ and $N = \max \{N_1, N_2\}$. Thus, the above inequalities (where $r_1$ and $r_2$ are swapped with $r$) hold simultaneously for $n\geq N$. Now, we have
\begin{align*}
    p_n\frac{a_n}{a_{2n}} - 2p_{2n} > r \quad &\text{and}\quad p_n\frac{a_n}{a_{2n+1}} - 2p_{2n+1} > r \\
    \Leftrightarrow  p_na_n - 2p_{2n}a_{2n} > ra_{2n} \quad &\text{and}\quad p_na_n - 2p_{2n+1}a_{2n+1} > ra_{2n+1}.
\end{align*}
Then for $M>N$, we have
\begin{align*}
    \sum_{n=N}^M(p_na_n - 2p_{2n}a_{2n}) > r\sum_{n=N}^Ma_{2n} \quad &\text{and}\quad \sum_{n=N}^M(p_na_n - 2p_{2n+1}a_{2n+1}) > r\sum_{n=N}^Ma_{2n+1}.
\end{align*}
Adding both inequalities together yields
\begin{align*}
    2\sum_{n=N}^Mp_na_n - &2\sum_{n=N}^M(p_{2n}a_{2n}+p_{2n+1}a_{2n+1})  > r\sum_{n=N}^M(a_{2n}+a_{2n+1}) \\
    \Leftrightarrow &2\sum_{n=N}^Mp_na_n - 2\sum_{n=2N}^{2M+1}p_{n}a_{n}  > r\sum_{n=2N}^{2M+1}a_{n} \\
    \Rightarrow &2\sum_{n=N}^{2M+1}p_na_n - 2\sum_{n=2N}^{2M+1}p_{n}a_{n}  > r\sum_{n=2N}^{2M+1}a_{n} \\
    \Rightarrow &2\sum_{n=N}^{2N-1}p_na_n > r\left(\sum_{n=1}^{2M+1}a_{n} - \sum_{n=1}^{2N-1}a_n\right) \\
    \Rightarrow &\frac{2\sum_{n=N}^{2N-1}p_na_n + r\sum_{n=N}^{2N-1}a_n}{r} > \sum_{n=1}^{2M+1}a_n.
\end{align*}
Since $N$ is a fixed number, then the left-hand side forms an upper bound on the partial sum. This implies that $\sum_{n=1}^\infty a_n$ converges. This proves \textit{(i)}.

Next we prove \textit{(ii)}. We only need to show that given the assumptions, then
\begin{align*}
    (\ref{eq3}) \Rightarrow \sum_{n=1}^\infty a_{2n}\;\;\text{diverges; or}
\end{align*}
\begin{align*}
    (\ref{eq4}) \Rightarrow \sum_{n=1}^\infty a_{2n+1}\;\;\text{diverges}.
\end{align*}
If at least one of these is satisfied, then we may conclude \textit{(ii)}. In particular, since
\begin{align*}
    \sum_{n=1}^\infty a_n = \sum_{n=1}^\infty a_{2n} + \sum_{n=0}^\infty a_{2n+1},
\end{align*}
we have $\sum_{n=1}^\infty a_{2n} \leq \sum_{n=1}^\infty a_n$. Thus, if (\ref{eq3}) holds, then $\sum_{n=1}^\infty a_n$ diverges by the ordinary comparison test. Similarly, since $\sum_{n=0}^\infty a_{2n+1} \leq \sum_{n=1}^\infty a_n$, then (\ref{eq4}) implies $\sum_{n=1}^\infty a_n$ diverges by the comparison test. Hence, we only need to prove the given conditional statements.

We will only prove $(\ref{eq3}) \Rightarrow \sum_{n=1}^\infty a_{2n}$ is divergent, since the proof for the second implication will be similar. Since we have $p_na_n< 2p_{2n}a_{2n}$ for sufficiently large $n$, say for all $n\geq N$ for some $N > 0$, then we define $k_0$ to be the smallest integer such that $2^{k_0} \geq N$. As a result, we have
\begin{align*}
    p_{2^k}a_{2^k}< 2p_{2^{k+1}}a_{2^{k+1}}, \quad \text{for all $k \geq k_0$}.
\end{align*}
If $k = k_0$, we have
\begin{align*}
    a_{2^{k_0+1}} > \frac{p_{2^{k_0}}}{2p_{2^{k_0+1}}}a_{2^{k_0}}.
\end{align*}
For $k = k_0+1$, we have
\begin{align*}
    a_{2^{k_0+2}} > \frac{p_{2^{k_0+1}}}{2p_{2^{k_0+2}}}a_{2^{k_0+1}} > \frac{p_{2^{k_0+1}}}{2p_{2^{k_0+2}}}\frac{p_{2^{k_0}}}{2p_{2^{k_0+1}}}a_{2^{k_0}} = \frac{p_{2^{k_0}}}{2^2p_{2^{k_0+2}}}a_{2^{k_0}}.
\end{align*}
Thus, we can sum up the series and use the inequality iteratively to obtain
\begin{align*}
    \sum_{k=k_0}^\infty a_{2^{k+1}} &\geq p_{2^{k_0}}a_{2^{k_0}}\left(\frac{1}{2p_{2^{k_0+1}}} + \frac{1}{2^2p_{2^{k_0+2}}} + ... \right) \\
    &= p_{2^{k_0}}a_{2^{k_0}}\sum_{k=k_0 +1}^\infty \frac{1}{2^{k-k_0}p_{2^k}} \\
    &= 2^{k_0}p_{2^{k_0}}a_{2^{k_0}}\sum_{k=k_0 +1}^\infty \frac{1}{2^{k}p_{2^k}}.
\end{align*}
By assumption, $\sum_{k=1}^\infty \frac{1}{2^{k}p_{2^k}}$ diverges. By the ordinary comparison test, we have $\sum_{k=k_0}^\infty a_{2^{k+1}}$ diverges. Since $\{a_{2^{k+1}}\}_{k=k_0}^\infty \subset \{a_{2n}\}_{n=1}^\infty$, we have $\sum_{n=1}^\infty a_{2n}$ diverges.
\end{proof}

\textbf{Remark.}
If \textit{(i)} holds for the sequence $p_n = 1$, we obtain the following corollary.

\textbf{Corollary 2.1.3}
Suppose $\{a_n\}$ is a sequence of positive numbers. If the following inequalities hold for positive numbers $r_1, r_2$ and $n$ sufficiently large:
\begin{align*}
    \frac{a_{2n}}{a_n} < \frac{1}{2+r_1}\quad \text{and}\quad \frac{a_{2n+1}}{a_n} < \frac{1}{2+r_2},
\end{align*}
then $\sum_{n=1}^\infty a_n$ converges.
Thus, we obtain a stricter form of the Second Ratio Test, since if the above conditions are met, then $L = \max \{ \limsup \frac{a_{2n}}{a_n}, \limsup \frac{a_{2n+1}}{a_n}\} < \frac{1}{2}$. 

On the other hand, letting $p_n = \frac{1}{n}$ in \textit{(ii)}, we obtain another corollary.

\textbf{Corollary 2.1.4}
Suppose $\{a_n\}$ is a sequence of positive numbers. If the either one of the following inequalities hold for $n$ sufficiently large:
\begin{align*}
    \frac{a_{2n}}{a_n} > 1\quad \text{or}\quad \frac{a_{2n+1}}{a_n} > 1 + \frac{1}{2n},
\end{align*}
then $\sum_{n=1}^\infty a_n$ diverges.

In comparison to the Second Ratio Test, it is sufficient to check whether one of the two ratios satisfies its corresponding inequality. In this way, this test is more general compared to the Second Ratio Test for divergence.

\section{Conclusion} 
The novelty of the Second Raabe's Test is that it addresses the inconclusive case found in the Second Ratio Test; nevertheless, like the Raabe's Test, there is also an inconclusive case. In this case, several tests were developed (the Second Gauss' Test, Second Kummer's Test) to partially address this case. However, more work could be done to further generalize the Second Gauss's Test and Second Kummer's Test, as we suspect that there exists a more general formulation. In particular, we speculate that there exists a characterization of convergence and divergence of series based on the Second Kummer's Test, as was shown in the ordinary Kummer's Test \cite{Tong1994}.

Ali \cite{Ali2008} constructed the $m$-th ratio test, which generalizes the Second Ratio Test from using two ratios to $m$ ratios, to test series convergence. While perhaps not too hard to formulate, another potential research direction would be the construction of series tests that address the inconclusive case found in the $m$-th ratio test.

\textbf{Acknowledgment.}
The author wishes to thank Professor Zhonghai Ding for his guidance on this work. The author would also like to thank Keoni Castellano, Michael Schwob, Shen Huang, and Bowen Liu for reading and proofreading the manuscript.


\vfill\eject


\begin{thebibliography}{1}
\bibitem{Ali2008} Ali, S. (2008). The \textit{m}th Ratio Test: New Convergence Tests for Series. \textit{Amer. Math. Monthly.}

\bibitem{d1768reflexions} d'Alembert, J. (1768). R{\'e}flexions sur les Suites Divergentes ou Convergentes. \textit{Opuscules Math{\'e}mathiques ou M{\'e}moires sur Diff{\'e}rens Sujets de G{\'e}om{\'e}trie.}

\bibitem{NotezurTheoriederConvergenzundDivergenzderReihen} Raabe, J. (1834). Note zur Theorie der Convergenz und Divergenz der Reihen. \textit{Journal für die reine und angewandte Mathematik.}

\bibitem{Tong1994} Tong, J. (1994). Kummer's Test Gives Characterization for Convergence or Divergence of all Positive Series. \textit{Amer. Math. Monthly.}
\end{thebibliography}
\end{document}